\documentclass[preprint,12pt]{elsarticle}

\usepackage{mathscinet}
\usepackage{amssymb}
\usepackage{amsmath}
\usepackage{amsthm}
\usepackage{esint}
\usepackage[colorlinks]{hyperref}
\AtBeginDocument{%
  \hypersetup{%
    linkcolor=blue,%
    citecolor=green,%
  }%
}
\usepackage{cleveref}
\usepackage[displaymath,mathlines,pagewise]{lineno}
\newcommand*\patchAmsMathEnvironmentForLineno[1]{%
  \expandafter\let\csname old#1\expandafter\endcsname\csname #1\endcsname
  \expandafter\let\csname oldend#1\expandafter\endcsname\csname end#1\endcsname
  \renewenvironment{#1}%
     {\linenomath\csname old#1\endcsname}%
     {\csname oldend#1\endcsname\endlinenomath}}%
\newcommand*\patchBothAmsMathEnvironmentsForLineno[1]{%
  \patchAmsMathEnvironmentForLineno{#1}%
  \patchAmsMathEnvironmentForLineno{#1*}}%
\AtBeginDocument{%
\patchBothAmsMathEnvironmentsForLineno{equation}%
\patchBothAmsMathEnvironmentsForLineno{align}%
\patchBothAmsMathEnvironmentsForLineno{flalign}%
\patchBothAmsMathEnvironmentsForLineno{alignat}%
\patchBothAmsMathEnvironmentsForLineno{gather}%
\patchBothAmsMathEnvironmentsForLineno{multline}%
}

\crefname{equation}{}{}

\newtheorem{theorem}{Theorem}[section]
\newtheorem{lemma}[theorem]{Lemma}

\newtheorem{corollary}[theorem]{Corollary}

\theoremstyle{definition}
\newtheorem{definition}[theorem]{Definition}

\newtheorem{notation}[theorem]{Notation}

\theoremstyle{remark}
\newtheorem{remark}[theorem]{Remark}

\numberwithin{equation}{section}

\journal{~}

\begin{document}

\begin{frontmatter}



\title{Interior pointwise $C^{\alpha}$ regularity for elliptic and parabolic equations with divergence-free drifts}

\author{Yuanyuan Lian}
\ead{lianyuanyuan.hthk@gmail.com; yuanyuanlian@correo.ugr.es}



\address[rvt]{Departamento de An\'{a}lisis Matem\'{a}tico,
Instituto de Matem\'{a}ticas IMAG, Universidad de Granada, Granada, Espa\~{n}a}

\begin{abstract}
We investigate the interior pointwise $C^{\alpha}$ regularity for weak solutions of elliptic and parabolic equations with divergence-free drifts. For such equations, the integrability condition on the drift can be relaxed and the interior $C^{\alpha}$ regularity for some $0<\alpha<1$ has been obtained previously with the aid of Harnack inequality. In this paper, we prove the interior pointwise $C^{\alpha}$ regularity for any $0<\alpha<1$ provided that the drift is small. We obtain the regularity under three different types conditions on the drift. The proof is based on the energy inequality and the perturbation technique.
\end{abstract}


\end{frontmatter}


\section{Introduction}
\label{intro}
This paper is devoted to investigate the interior pointwise regularity for weak solutions of elliptic and parabolic equations:
\begin{equation}\label{e.equation}
  \begin{aligned}
    &-\Delta u+b^iu_i=0 ~~\mbox{in}~~B_1,&&~~\mbox{(\textbf{E})},\\
    &u_t-\Delta u+b^iu_i=0 ~~\mbox{in}~~Q_1,&&~~\mbox{(\textbf{P})}.~~
  \end{aligned}
\end{equation}
Here, $B_1\subset \mathbb{R}^n$ is the unit ball and $Q_1\subset \mathbb{R}^{n+1}$ is the unit parabolic cube ($n\geq 2$). In above equations, the Einstein summation convention is used (similarly hereinafter).

In this paper, we consider weak solutions $u\in H^1(B_1)$ and $u\in L^2(-1,0;H^1(B_1))\cap L^{\infty}(-1,0;L^2(B_1))$ respectively in (\textbf{E}) and (\textbf{P}). Additionally, we always assume the following divergence-free structure:
\begin{equation}\label{e.divb0}
\mathrm{div}~\pmb{b} = 0~~\mbox{in weak sense.}
\end{equation}
That is,
\begin{equation*}
  \begin{aligned}
    &\int_{B_1} b^i \varphi_i =0,~~\forall \varphi \in C^{\infty}_c(B_1)~~\mbox{for (\textbf{E})},\\
    &\int_{Q_1} b^i \varphi_i =0,~~\forall \varphi \in C^{\infty}_c(Q_1)~~\mbox{for (\textbf{P})}.
  \end{aligned}
\end{equation*}

The H\"{o}lder regularity for weak solutions of elliptic and parabolic equations in divergence form have been studied extensively from the beginning of De Giorgi \cite{MR0093649}, Nash \cite{MR0100158} and Moser \cite{MR0159138, MR0206469}. For equations with the lower order term $b^iu_i$, Morrey \cite{MR120446} got the H\"{o}lder regularity in the elliptic case if $\pmb{b}\in L^p$ $(p\geq n)$; Ladyzhenskaya and Ural'tseva \cite{MR141891} obtained the corresponding result in the parabolic case if $\pmb{b}\in L^{p+2} (p>n)$. Stampacchia \cite{MR192177} also got the H\"{o}lder continuity with $\pmb{b} \in L^n$ for the elliptic equation. Trudinger \cite{MR0369884} proved Harnack inequality in the elliptic case with $\pmb{b}\in L^p (p>n)$ and Lieberman \cite[Chapter 6]{MR1465184} extended the result in the parabolic case with $\pmb{b}\in L^{p,p'} (n/p+2/p'<1)$. It is well known that the Harnack inequality implies the H\"{o}lder regularity.

Equations with lower order coefficients satisfying \Cref{e.divb0} arise in some problems such as incompressible flows (see \cite{MR2753290, MR2660718, MR2545826, MR1736969, MR2852216}). Researchers have shown particular interest in the quantitative properties of solutions for elliptic and parabolic equations with divergence-free drifts (see \cite{MR2784072, MR3164858, MR3462128, MR2760150, MR2031029}). Ignatova, Kukavica and Ryzhik \cite{MR3164858, MR3462128} proved the Harnack inequality for elliptic and parabolic equations under $b\in L^{p}~(n/2<p<n)$ and $b\in L^{p,p'}~(1<n/p+2/p'<2)$ respectively. However, the constant in this Harnack inequality is not universal and the solution is a Lipschitz solution rather than a weak solution. Nazarov and Ural'tseva \cite{MR2760150} obtained the H\"{o}lder regularity for some $0<\alpha<1$ under the assumptions:
\begin{equation}\label{e1.2}
 \begin{aligned}
 &\pmb{b}\in L^{p}(B_1)~~(\frac{n}{2}<p<n, p\geq 2),
 ~~\pmb{b}\in C_p^{-1}(B_1),~~&&(\mbox{E}),\\
 &\pmb{b}\in L^{p,p'}(Q_1)~~(1<\frac{n}{p}+\frac{2}{p'}<2),
 ~~\pmb{b}\in C_{p,p'}^{-1}(Q_1),~~&&(\mbox{P}).
 \end{aligned}
\end{equation}
The (E) (resp. (P)) is for the elliptic equations (resp. parabolic equations) (similarly hereinafter). The $C_p^{-1}(B_1)$ (resp. $C_{p,p'}^{-1}(Q_1)$) (see \Cref{d-f2}) is expressed by $\pmb{b}\in \mathbb{M}^{n/p-1}_p(B_1)$ (resp. $\pmb{b}\in \mathbb{M}^{n/p+2/p'-1}_{p,p'}(Q_1)$) in \cite{MR2760150}, which is a scale of Morrey spaces. Liskevich and Zhang \cite{MR2128546} and Zhang \cite{MR2031029} obtained the H\"{o}lder continuity to the elliptic and parabolic equations for some $0<\alpha<1$ with the following assumptions
\begin{equation}\label{e1.3}
\begin{aligned}
 &\pmb{b}\in L^2(B_1),~~\int_{B_1} |\pmb{b}|^2 \varphi^2 \leq C_{E} \int_{B_1} |D \varphi|^2,
      ~~\forall~\varphi\in C^{\infty}_c (B_1),~~&&(\mbox{E}),\\
 &\pmb{b}\in L^2(Q_1),~~\int_{Q_1} |\pmb{b}|^2 \varphi^2 \leq C_{P} \int_{Q_1} |D \varphi|^2,
      ~~\forall~\varphi\in C^{\infty}_c (Q_1),~~&&(\mbox{P}).
\end{aligned}
\end{equation}
In fact, Zhang \cite{MR2031029} assumed $\pmb{b}=\pmb{b}(x)$ in (P) of \cref{e1.3}. Seregin et al. \cite{MR2852216} proved the H\"{o}lder regularity for some $0<\alpha<1$ if $\pmb{b}\in BMO^{-1}(B_1)$ (resp. $\pmb{b}\in L^{\infty}(-1,0;BMO^{-1}(B_1))$) for elliptic equations (resp. parabolic equations). It means that $b^i=(A^{ij})_j$ where
\begin{equation}\label{e1.4}
  \begin{aligned}
    &A^{ij}=-A^{ji},~ \pmb{A}\in BMO(B_1),~~&&(\mbox{E}),\\
    &A^{ij}=-A^{ji},~\pmb{A}\in L^{\infty}(-1,0;BMO(B_1)),~~&&(\mbox{P}).
  \end{aligned}
\end{equation}
Li and Pipher \cite{MR3936348} provided a simpler proof of the interior $C^{\alpha}$ estimate for the elliptic case.  Friedlander and Vicol \cite{MR2784072} proved the H\"{o}lder continuity of weak solutions to parabolic equations with $\pmb{b}\in L^2(Q_1)\cap L^{\infty}(-1,0;BMO^{-1}(B_1))$.

One key in all above results is that with the aid of \Cref{e.divb0}, the integrability condition on $\pmb{b}$ can be relaxed when multiplying test functions. We explain this formally. Take the elliptic equation for example. By taking the test function $\varphi^2 u$ where $\varphi \in C^{\infty}_c(B_1)$ in the definition of weak solutions, the lower term $b^iu_i$ can be treated as follows:
\begin{equation*}
  -\int b^i u_i \varphi^2 u=\int b^i_i u^2 \varphi^2+b^i u_i (2 \varphi \varphi_i) u+b^i u \varphi^2 u_i.
\end{equation*}
By $\mathrm{div}~\pmb{b}=0$,
\begin{equation*}
-\int b^i u_i \varphi^2 u=\int b^i\varphi \varphi_i u^2.
\end{equation*}
Since $u\in H^{1}$, to ensure the left-hand term making sense, we must require $\pmb{b}\in L^n$. However, $\pmb{b}\in L^p~(p>n/2)$ is enough to ensure the right-hand term making sense. Certainly, the above formula is just a formal calculation but it reveals the essence how $\mathrm{div}~\pmb{b}= 0$ relaxes the integrability condition on $\pmb{b}$.


The H\"{o}lder continuity can be regarded equivalently as scaling invariant decay. That is, $u\in C^{\alpha}(0)$ is equivalent to
\begin{equation*}
  \underset{B_{\eta r}}{\mathrm{osc}}~u\leq \mu \underset{B_r}{\mathrm{osc}}~u,~~\forall~0<r<1,
\end{equation*}
where $0<\eta, \mu<1$ are two constants independent of $r$. Hence, to obtain the H\"{o}lder regularity, we must require that the equation is scaling invariant. Precisely, consider the transformation $y=x/r, v(y)=u(x)/r^{\alpha}$. Then $v$ is a weak solution of
\begin{equation*}
  -\Delta v+\tilde{b}^iv_i=0,
\end{equation*}
where $\tilde{\pmb{b}}(y)=r\pmb{b}(x)$. The scaling invariance is that $\tilde{\pmb{b}}$ is the same as (or better than) $\pmb{b}$ in some sense. A typical example is the $L^n$ norm, i.e., $\|\tilde{\pmb{b}}\|_{L^n}=\|\pmb{b}\|_{L^n}$. From the viewpoint of scaling, $\pmb{b}$ has the scaling $x^{-1}$, or $\pmb{b}\in C^{-1}$. In fact, \cref{e1.2}, \cref{e1.3} and \cref{e1.4} are three different types of characterization of $\pmb{b}\in C^{-1}$.

All regularity results obtained in above literatures under the condition \cref{e1.2}, \cref{e1.3} or \cref{e1.4} are $C^{\alpha}$ for some $0<\alpha<1$ since the proofs are based on the Harnack inequality. In this paper, we prove the interior pointwise $C^{\alpha}$ regularity for \textbf{any} $0<\alpha<1$ by enhancing that $\pmb{b}$ is ``small'', i.e.,
\begin{equation}\label{e1.5}
 \begin{aligned}
 &\|\pmb{b}\|_{C_p^{-1}(0)}<\delta,~~&&\|\pmb{b}\|_{C_{p,p'}^{-1}(0)}<\delta~~\mbox{in}~~\Cref{e1.2},\\
 &C_E<\delta,~~&&C_P<\delta~~\mbox{in}~~\Cref{e1.3},\\
 &\|\pmb{A}\|_{BMO_q(0)}<\delta,~~&&\|\pmb{A}\|_{L^{\infty}(-1,0;BMO_q(0))}<\delta~~\mbox{in}~~\Cref{e1.4},
 \end{aligned}
\end{equation}
where $\delta$ depends only on $n,\alpha$ (also $p$ in \Cref{e1.2} or $q$ in \Cref{e1.4}). Here $BMO_q(0)$ (see \Cref{d-f3}) is a pointwise BMO condition which is enough to prove the pointwise $C^{\alpha}$ regularity. With the ``small'' $\pmb{b}$, (\textbf{E}) and (\textbf{P}) can be considered as perturbations of the Laplace's equation and the heat equation respectively. Since the energy inequalities can provide the compactness, the H\"{o}lder regularity can be obtained by the compactness method and the scaling invariance argument.
~\\

Now, we introduce some notions and definitions. Let $s,s'\geq1$, $\Omega \subset \mathbb{R}^{n}$ be a bounded domain and $U:=\Omega\times (T_1,T_2)$. Define
\begin{equation*}
\|f\|^*_{L^{s}(\Omega)}:=\left(\fint_{\Omega} |f|^s\right)^{1/s}
~~\mbox{and}~~
\|f\|^*_{L^{s,s'}(U)}:=\left(\fint^{T_2}_{T_1}  \left(\fint_{\Omega} |f|^s\right)^{s'/s} \right)^{1/s'}
\end{equation*}
for elliptic and parabolic equations respectively, where $\fint$ denotes the integral average. The benefit of this notion is that the scaling of $\|f\|^*_{L^{s}(\Omega)}$ and $\|f\|^*_{L^{s,s'}(U)}$ are exactly that of $f$, which can make the proof concise when we deal with the scaling argument.

\begin{remark}\label{re1.1}
The $f\in L^{s,s'}(U)$ means
\begin{equation*}
\|f\|_{L^{s,s'}(U)}:=\left(\int_{T_1}^{T_2}\|f(t)\|^{s'}_{L^s(\Omega)}dt\right)^{\frac{1}{s'}}
:=\left(\int_{T_1}^{T_2}\left(\int_{\Omega}|f(x,t)|^sdx\right)^{\frac{s'}{s}}dt\right)^{\frac{1}{s'}}<\infty,
\end{equation*}
which is used widely for parabolic equations (see \cite{MR0241822,MR0150444}). It is due to the occurrence of a function having different integrability with respect to $x$ and $t$ when treating parabolic equations (e.g., the weak solution $u$).
\end{remark}
~\\

Next, we introduce some definitions of pointwise properties for a function.
\begin{definition}\label{d-f1}
Let $s\geq 1$ and $f:\Omega \rightarrow \mathbb{R}$ be a function and $0<\alpha\leq 1$. We say that $f$ is $C_{s}^{\alpha}$ at $x_0\in \Omega$ or $f\in C_{s}^{\alpha}(x_0)$ if there exist $K,r_0>0$ and a constant $P$ such that
\begin{equation}\label{m-holder}
  \|f-P\|^*_{L^{s}(B_{r}(x_0)\cap \Omega)}\leq K r^{\alpha},~~\forall~0<r<r_0.
\end{equation}
Then we define
\begin{equation*}
[f]_{C_{s}^{\alpha}(x_0)}=\min \left\{K\big | \cref{m-holder} ~\mbox{holds with}~K\right\}
\end{equation*}
and
\begin{equation*}
\|f\|_{C_{s}^{\alpha}(x_0)}=|P|+[f]_{C_{s}^{\alpha}(x_0)}.
\end{equation*}
If $f\in C_{s}^{\alpha}(x)$ for any $x \in \Omega'\subset \Omega$ with the same $r_0$ and
\begin{equation*}
  \|f\|_{C_{s}^{\alpha}(\bar{\Omega}')}:= \sup_{x \in \Omega'} \|f\|_{C_{s}^{\alpha}(x)}<+\infty,
\end{equation*}
we say that $f\in C_{s}^{\alpha}(\bar{\Omega}')$.
\end{definition}

\begin{remark}\label{r-df1.2}
It is well known if $\Omega$ is a Lipschitz domain, the above definition $f\in C_{s}^{\alpha}(\bar{\Omega}')$ is equivalent to the classical definition of $f\in C^{\alpha}(\bar{\Omega})$.
\end{remark}
~\\

Additionally, we define the following type of continuity.
\begin{definition}\label{d-f2}
Let $1\leq s\leq n$. We say that $f$ is $C_{s}^{-1}$ at $x_0\in \Omega$ or $f\in C_{s}^{-1}(x_0)$ if there exist $K,r_0>0$ such that
\begin{equation}\label{e.c-1}
\|f\|^*_{L^{s}(B_r(x_0)\cap \Omega)}\leq K r^{-1}, ~\forall ~0<r<r_0.
\end{equation}
Then define
\begin{equation*}
\|f\|_{C_{s}^{-1}(x_0)}=\min \left\{K \big | \cref{e.c-1} ~\mbox{holds with}~K\right\}.
\end{equation*}
If $f\in C_{s}^{-1}(x)$ for any $x\in \Omega'\subset \Omega$ with the same $r_0$ and
\begin{equation*}
  \|f\|_{C_{s}^{-1}(\bar{\Omega}')}:= \sup_{x\in \Omega'} \|f\|_{C_{s}^{-1}(x)}<+\infty,
\end{equation*}
we say that $f\in C_{s}^{-1}(\bar{\Omega}')$.
\end{definition}

\begin{remark}\label{r-df2.2}
In this paper, if we use \Cref{d-f1} and \Cref{d-f2} for some functions, we always assume that $r_0=1$ without loss of generality.
\end{remark}

\begin{remark}\label{r-df2.21}
If $\pmb{b}\in L^n$, then $\pmb{b}\in C^{-1}_s$ for any $1\leq s \leq n$.
\end{remark}

\begin{remark}\label{r-df2.3}
If $\pmb{f}$ is vector valued function, $\pmb{f}\in C_{s}^{-1}(x_0)$ means that each component of $\pmb{f}$ is $C_{s}^{-1}$ at $x_0$ and $\|\pmb{f}\|_{C_{s}^{-1}(x_0)}$ can be defined correspondingly.
\end{remark}

\begin{remark}\label{rdf2.4}
For parabolic equations, we can define the pointwise smoothness correspondingly and we only need to replace $\|\cdot\|^*_{L^{s}(B_{r}(x_0))}$ by $\|\cdot\|^*_{L^{s,s'}(Q_{r}(x_0,t_0))}$ in \cref{m-holder} and \cref{e.c-1}.
\end{remark}

\begin{remark}\label{r-df2.5}
Since $\pmb{b}\in L^{p}(B_1)$ (resp. $\pmb{b} \in L^{p,p'}(Q_1)$) where $p$ (resp. $p,p'$) is fixed in this paper, the subscripts are omitted when we state the pointwise properties (in \Cref{d-f1} and \Cref{d-f2}) of $\pmb{b}$, for example, we write $\pmb{b}\in C^{-1}(0)$ in place of $\pmb{b}\in C^{-1}_{p}(0)$ etc.
\end{remark}
~\\

Furthermore, we define the pointwise $BMO$ regularity.
\begin{definition}\label{d-f3}
We say that $f$ is $BMO_s$ at $x_0\in \Omega$ or $f\in BMO_s(x_0)$ with radius $r_0>0$ if
\begin{equation*}
   \|f\|_{BMO_s(x_0)}:=\sup_{0<r<r_0} \| f-f_{B_r(x_0)\cap \Omega}\|^*_{L^s(B_r(x_0)\cap \Omega)}<+\infty,
\end{equation*}
where $f_{B_r(x_0)\cap \Omega}=\fint_{B_r(x_0)\cap \Omega} f$. If $f\in BMO_s(x_0)$ for any $x_0\in \Omega$ with the same radius $r_0$ and
\begin{equation*}
 \|f\|_{BMO_s(\Omega)}:=\sup_{x\in \Omega}\|f\|_{BMO(x)} <+\infty,
\end{equation*}
we say that $f\in BMO_s(\Omega)$.

Let $g=\mathrm{div}~\pmb{f}$, we say that $g\in BMO_s^{-1}(0)$ if $\pmb{f}\in BMO_s(0)$ and define
\begin{equation*}
  \|g\|_{BMO_s^{-1}(0)}=\|\pmb{f}\|_{BMO_s(0)}.
\end{equation*}
Moreover, we define $g\in BMO_s^{-1}(\Omega)$ if $\pmb{f}\in BMO_s(\Omega)$ and
\begin{equation*}
  \|g\|_{BMO_s^{-1}(\Omega)}=\|\pmb{f}\|_{BMO_s(\Omega)}.
\end{equation*}
\end{definition}

\begin{remark}\label{r-df3.1}
Note that the usual BMO space is defined based on $L^1$ norm rather than $L^s$ norm. In fact, these definitions are equivalent (see \cite[Corollary on P. 144]{MR1232192} and \cite[Corollary 6.22]{MR3099262}).
\end{remark}

\begin{remark}\label{r-df3.2}
For a vector valued function $\pmb{g}$, we say $\pmb{g}\in BMO^{-1}$ if each component $g^i\in BMO^{-1}$.
\end{remark}

\begin{remark}
  Since we assume that $\pmb{A}\in BMO_q$ ($q>n$) where $q$ is fixed in this paper, the subscripts are omitted when we state the BMO properties, i.e., $\pmb{A}\in BMO_q(0)$ is abbreviated as $\pmb{A}\in BMO(0)$.
\end{remark}
~\\

Next, we give definitions for weak solutions.
%

\begin{definition}\label{Weak-ZQ}
Let $\pmb{b}\in L^2(B_1)$. We say that $u\in H^1(B_1)$ is a weak solution of (\textbf{E}) if
\begin{equation}\label{e.weak1}
\int_{B_1}u_i \varphi_i+b^iu_i \varphi=0,~~\forall~\varphi\in C_c^{\infty}(B_1).
\end{equation}
\end{definition}

\begin{remark}\label{re.we.u}
If $\pmb{b}\in L^{\infty}(B_1)$, $\varphi$ can be chosen in $H^1_0(B_1)$.
\end{remark}
~\\

\begin{definition}\label{Weak_BMO}
With $\pmb{b}\in BMO^{-1}(B_1)$, we say that $u\in H^1(B_1)$ is a weak solution of (\textbf{E}) if
\begin{equation}\label{e.weak2}
\int_{B_1}u_i \varphi_i-A^{ij}u_i \varphi_j=0,~~\forall~\varphi\in H_0^1(B_1).
\end{equation}
where $A^{ij}u_i\varphi_j :=  \frac{1}{2} A^{ij}(u_i \varphi_j-u_j \varphi_i)$.
\end{definition}

\begin{remark}\label{re.W.2}
Since $A^{ij}\in BMO$ and $u_i \varphi_j-u_j \varphi_i \in \mathcal{H}^1$ (the Hardy space), the above definition makes sense (see \cite{MR1225511} and \cite[P. 160]{MR3936348}).
\end{remark}

\begin{remark}\label{re.W.3}
We can use
\begin{equation}\label{e.weak2.1}
\int_{B_1}u_i \varphi_i-(A^{ij}-A^{ij}_{B_1})u_i \varphi_j=0
\end{equation}
to replace \cref{e.weak2}. Since $A^{ij}=-A^{ji}$, we have $\int_{B_1} A^{ij}_{B_1}u_i \varphi_j=0$.
\end{remark}

\begin{remark}\label{re-W.1}
If $\pmb{b}\in L^{p,p'}$, we refer to \cite[P. 136]{MR0241822} for the definition of weak solutions with regard to parabolic equations. The weak solutions corresponding to $\pmb{b}\in L^{\infty}(-1,0;BMO^{-1}(B_1))$ for the parabolic case can be defined similarly to \cref{e.weak2}.
\end{remark}
~\\

Our main results are in the following. Since we consider the pointwise regularity, we focus on the regularity at the origin $0$ throughout this paper. We use the $L^2$ norm in \Cref{m-holder} for the solution $u$ and we write $u\in C^{\alpha}$ instead of $u\in C_2^{\alpha}$ for simplicity.

First, we give the $C^{\alpha}$ regularity for the elliptic equations.
\begin{theorem}\label{t-Ca-elli}
Let $0<\alpha<1$ and $u\in H^1(B_1)$ be a weak solution of (\textbf{E}) in \Cref{e.equation}. Suppose that
\begin{equation*}
  \begin{aligned}
    &\|\pmb{b}\|_{C^{-1}(0)}\leq \delta~~(\mbox{or}~~C_E\leq \delta
   ~~\mbox{or}~~\|\pmb{b}\|_{BMO^{-1}(0)}\leq \delta),
  \end{aligned}
\end{equation*}
where $\delta$ is small and depends only on $n, p,\alpha$ (or $n, \alpha$ or $n, q,\alpha$) (see \Cref{e1.5}).
Then $u\in C^{\alpha}(0)$, i.e., there exists a constant $P$ such that
\begin{equation}\label{e.Ca.esti}
  \|u-P\|^*_{L^2(B_r)}\leq C r^{\alpha}\|u\|^*_{L^{2}(B_1)}, ~~\forall ~0<r<1
\end{equation}
and
\begin{equation*}
  |P|\leq C\|u\|^*_{L^{2}(B_1)},
\end{equation*}
where $C$ depends only on $n, p,\alpha$ (or $n, \alpha$ or $n, q,\alpha$).
\end{theorem}

As a direct consequence, we have the following corollary.
\begin{corollary}\label{coro.1}
Let $0<\alpha<1$ and $u\in H^1(B_1)$ be a weak solution of (\textbf{E}) in \Cref{e.equation}. Suppose that
\begin{equation*}
  \begin{aligned}
    &\|\pmb{b}\|_{C^{-1}(\bar B_1)}\leq \delta~~(\mbox{or}~~C_E\leq \delta
   ~~\mbox{or}~~\|\pmb{b}\|_{BMO^{-1}(\bar B_1)}\leq \delta),
  \end{aligned}
\end{equation*}
where $\delta$ is small and depends only on $n, p,\alpha$ (or $n, \alpha$ or $n, q,\alpha$).
Then $u\in C^{\alpha}(\bar B_{1/2})$ and
\begin{equation}\label{e.Ca.esti}
  \|u\|_{C^{\alpha}(\bar B_{1/2})}=\sup_{x,y\in B_{1/2},x\neq y} \frac{|u(x)-u(y)|}{|x-y|^{\alpha}}\leq C \|u\|_{L^{2}(B_1)},
\end{equation}
where $C$ depends only on $n, p,\alpha$ (or $n, \alpha$ or $n, q,\alpha$).
\end{corollary}

Next, the $C^{\alpha}$ regularity for the parabolic equation is:
\begin{theorem}\label{t-Ca-para}
Let $0<\alpha<1$ and $u\in L^2(-1,0;H^1(B_1))\cap L^{\infty}(-1,0;L^2(B_1))$ be a weak solution of (\textbf{P}) in \Cref{e.equation}. Suppose that
\begin{equation*}
  \begin{aligned}
    &\|\pmb{b}\|_{C^{-1}(0)}\leq \delta~~(\mbox{or}~~C_P\leq \delta
   ~~\mbox{or}~~\|\pmb{b}\|_{L^{\infty}(-1,0;BMO^{-1}(0))}\leq \delta)
  \end{aligned}
\end{equation*}
where $\delta$ is small and depends only on $n, p, \alpha$ (or $n, \alpha$ or $n, q, \alpha$).
Then $u\in C^{\alpha}(0)$, i.e., there exists a constant $P$ such that
\begin{equation}\label{e.Ca.esti}
  \|u-P\|^*_{L^2(Q_r)}\leq C r^{\alpha}\|u\|^*_{L^{2}(Q_1)}, ~~\forall ~0<r<1
\end{equation}
and
\begin{equation*}
  |P|\leq C\|u\|^*_{L^{2}(Q_1)},
\end{equation*}
where $C$ depends only on $n, p,\alpha$ (or $n, \alpha$ or $n, q,\alpha$).
\end{theorem}

Similarly, for the parabolic equation, we can obtain the following corollary.
\begin{corollary}\label{coro.2}
Let $0<\alpha<1$ and $u\in L^2(-1,0;H^1(B_1))\cap L^{\infty}(-1,0;L^2(B_1))$ be a weak solution of (\textbf{P}) in \Cref{e.equation}. Suppose that
\begin{equation*}
  \begin{aligned}
    &\|\pmb{b}\|_{C^{-1}(\bar{Q}_1)}\leq \delta~~(\mbox{or}~~C_P\leq \delta
   ~~\mbox{or}~~\|\pmb{b}\|_{L^{\infty}(-1,0;BMO^{-1}(\bar{B}_1))}\leq \delta)
  \end{aligned}
\end{equation*}
where $\delta$ is small and depends only on $n, p, \alpha$ (or $n, \alpha$ or $n, q, \alpha$).
Then $u\in C^{\alpha}(\bar{Q}_{1/2})$ and
\begin{equation}\label{e.Ca.esti}
  \|u\|_{C^{\alpha}(\bar{Q}_{1/2})}\leq C \|u\|_{L^{2}(Q_1)},
\end{equation}
where $C$ depends only on $n, p,\alpha$ (or $n, \alpha$ or $n, q,\alpha$).
\end{corollary}
~\\

The paper is organized as follows. In \Cref{sec:1}, we first prove the energy inequality for elliptic equations under conditions \cref{e1.2}, \cref{e1.3} or \cref{e1.4}. Then we obtain the interior pointwise $C^{\alpha}$ regularity by the compactness method and a scaling argument. For the parabolic equations, the interior pointwise $C^{\alpha}$ regularity is obtained in \Cref{sec:2}.~\\

\begin{notation}\label{no1.1}
\begin{enumerate}~~\\
\item $|x|:=\left(\sum_{i=1}^{n} x_i^2\right)^{1/2}$ for $x=(x_1,x_2,...,x_n) \in \mathbb{R}^n$  and $|(x,t)|:=(|x|^2+|t|)^{1/2}$ for $(x,t)\in \mathbb{R}^{n+1}$.
\item $B_r(x_0)=B(x_0,r)=\{x\in \mathbb{R}^{n}\big| |x-x_0|<r\}$ and $B_r=B_r(0)$.
\item $Q_r(x_0,t_0)=Q((x_0,t_0),r)=B_r(x_0)\times (t_0-r^2,t_0]$ and $Q_r=Q_r(0)$.
\item $\bar A $: the closure of $A$, where $A\subset \mathbb{R}^{n}$ or $\mathbb{R}^{n+1}$.
\item $\varphi_t=D_t \varphi$, $\varphi _i=D_i \varphi=\partial \varphi/\partial x _{i}$ $(1\leq i\leq n)$, $\varphi _{ij}=D_{ij}\varphi =\partial ^{2}\varphi/\partial x_{i}\partial x_{j}$ $(1\leq i,j\leq n)$.
\end{enumerate}
\end{notation}

\section{$C^{\alpha}$ regularity for elliptic equations}\label{sec:1}
In this section, we give the proof of the main result for the elliptic equation (\textbf{E}) in \Cref{e.equation}. For brevity, we just say that $u$ is a weak solution of \Cref{e.equation} in this section. The compactness method is applied in this paper. Hence, we first prove the energy inequality (Caccioppoli inequality), which provides the necessary compactness.

\begin{lemma}\label{ener-e}
  Let $u$ be a weak solution of \Cref{e.equation}. Assume that $\pmb{b}$ satisfies one of the following conditions:
  \begin{equation*}
    \begin{aligned}
      &(1)~\pmb{b}\in L^{p}(B_1)~~(\frac{n}{2}<p<n, p\geq 2);\\
      &(2)~\pmb{b}\in L^2(B_1),~~\int_{B_1} |\pmb{b}|^2 \varphi^2 \leq C_{E} \int_{B_1} |D \varphi|^2,
      ~~\forall~\varphi\in C^{\infty}_c (B_1);\\
      &(3)~\pmb{b}\in BMO^{-1}(B_1).
    \end{aligned}
  \end{equation*}
  Then
  \begin{equation}\label{e1.6}
    \|Du\|_{L^2(B_{1/2})}\leq C \|u\|_{L^2(B_{1})},
  \end{equation}
  where $C$ depends only on $n$ and $\|\pmb{b}\|_{L^p(B_1)}$ (or $C_E$ or $\|\pmb{b}\|_{BMO^{-1}(0)}$).
\end{lemma}

\proof
\textbf{Case (1):} $\pmb{b}\in L^{p}(B_1)$ ($n/2<p<n$, $p\geq 2$). We first assume that $\pmb{b}\in C^{\infty}(\bar{B}_1)$. Let
\begin{equation*}
\psi=\left(1-|x|^2\right)^{\beta},
\end{equation*}
where $\beta$ is large to be specified later. Take $\psi^2 u$ as test function in \cref{e.weak1} to get
\begin{equation*}
  \int_{B_1} \psi^2 u_i^2+2\psi \psi_i u u_i+b^i u_i \psi^2 u=0.
\end{equation*}
With the aid of the divergence-free structure condition \cref{e.divb0},
\begin{equation}\label{e.energy.1}
  \int_{B_1} \psi^2 u_i^2=-\int_{B_1}2\psi \psi_i u u_i+b^i u_i \psi^2 u =-2\int_{B_1}\psi \psi_i u u_i + \int_{B_1} b^i u^2 \psi \psi_i.
\end{equation}

In the following, we estimate each term of the right hand. By the Young’s inequality,
\begin{equation*}
  \begin{aligned}
    \left|\int_{B_1} 2\psi \psi_i u u_i\right|
    &\leq \frac{1}{8} \int_{B_1} \psi^2 u_i^2+\bar C \int_{B_1} \psi_i^2 u^2.
  \end{aligned}
\end{equation*}
Choose $0<\delta_0\leq \frac{4}{n+2}$ and $\beta\geq \frac{2}{\delta_0}$. Then we have
\begin{equation*}
\frac{|D \psi|}{\psi^{1-\frac{\delta_0}{2}}}\leq C ~~\mbox{ in}~~B_1.
\end{equation*}
By the H\"{o}lder inequality, the Sobolev inequality and the Young’s inequality,
\begin{equation*}
  \begin{aligned}
    \int_{B_1} \psi_i^2 u^2
    &=\int_{B_1} \psi^{2-\delta_0} |u|^{2-\delta_0} \frac{\psi_i^2}{\psi^{2-\delta_0}} |u|^{\delta_0}\\
    &\leq C\left(\int_{B_1} \left(\frac{|D\psi|}{\psi^{1-\frac{\delta_0}{2}}}\right)^{\frac{2}{\delta_0}} |u| \right)^{\delta_0}
    \left(\int_{B_1}(\psi |u|)^{\frac{2-\delta_0}{1-\delta_0}}\right)^{1-\delta_0}\\
    &\leq C\|u\|_{L^1(B_1)}^{\delta_0}
    \left(\|\psi Du\|_{L^2(B_1)}^{2-\delta_0}+\|u D\psi \|_{L^2(B_1)}^{2-\delta_0}\right)\\
    &\leq \varepsilon \|\psi Du\|_{L^2(B_1)}^{2}+\varepsilon \|u D\psi\|_{L^2(B_1)}^{2}
    +C \|u\|_{L^1(B_1)}^{2},
  \end{aligned}
\end{equation*}
where $\varepsilon$ is small enough such that
\begin{equation*}
  \frac{\varepsilon \bar{C}}{1-\varepsilon}\leq \frac{1}{8}.
\end{equation*}
Then
\begin{equation}\label{e2.11}
  \begin{aligned}
   \bar{C} \int_{B_1} \psi_i^2 u^2
    \leq \frac{1}{8} \|\psi Du\|_{L^2(B_1)}^{2}+C \|u\|_{L^1(B_1)}^{2}.
  \end{aligned}
\end{equation}
Thus,
\begin{equation*}
  \begin{aligned}
    \left|\int_{B_1} 2\psi \psi_i u u_i\right|
    &\leq \frac{1}{4} \|\psi Du\|_{L^2(B_1)}^{2}+C \|u\|_{L^1(B_1)}^{2}.
  \end{aligned}
\end{equation*}

Next, take $0<\delta_0<1$ such that (note that $p>n/2$)
\begin{equation*}
  \frac{p(2-\delta_0)}{p(1-\delta_0)-1}\leq \frac{2n}{n-2}.
\end{equation*}
Then by setting $\beta\geq \frac{2}{\delta_0}$,
\begin{equation*}
\frac{|D \psi|}{\psi^{1-\delta_0}}\leq C ~~\mbox{ in}~~B_1.
\end{equation*}
Hence,
\begin{equation}\label{e.1.b}
  \begin{aligned}
    &\left|\int_{B_1}b^i u^2 \psi \psi_i\right|
    =\left|\int_{B_1}b^i |u|^{2-\delta_0} \psi^{2-\delta_0} \frac{\psi_i}{\psi^{1-\delta_0}}|u|^{\delta_0}\right|\\
    &\leq C  \left(\int_{B_1}\left(\frac{|D \psi|}{\psi^{1-\delta_0}}\right)^{\frac{1}{\delta_0}} |u| \right)^{\delta_0} \left(\int_{B_1}(|\pmb{b}| |u|^{2-\delta_0} \psi^{2-\delta_0})^{\frac{1}{1-\delta_0}}\right)^{1-\delta_0}   \\
    &\leq C\|u\|_{L^1(B_1)}^{\delta_0}
    \|\pmb{b}\|_{L^p(B_1)}
     \left(\int_{B_1} \left(u \psi\right)^{\frac{p(2-\delta_0)}{p(1-\delta_0)-1}}\right)^{\frac{p(1-\delta_0)-1}{p}}\\
    &\leq C\|\pmb{b}\|_{L^p(B_1)} \|u\|_{L^1(B_1)}^{\delta_0}
    \left(\|\psi Du\|_{L^2(B_1)}^{2-\delta_0}+\|u D\psi\|_{L^2(B_1)}^{2-\delta_0}\right)\\
    &\leq \frac{1}{8} \|\psi Du\|_{L^2(B_1)}^{2}+\frac{1}{8}\|u D\psi\|_{L^2(B_1)}^{2}+C \|u\|_{L^1(B_1)}^{2}\\
    &\leq \frac{1}{4} \|\psi Du\|_{L^2(B_1)}^{2}+C \|u\|_{L^1(B_1)}^{2},
  \end{aligned}
\end{equation}
where $C$ depends only on $n,p$ and $\|\pmb{b}\|_{L^p(B_1)}$.
Therefore, by combing the above estimates, we have
\begin{equation*}\label{e.vk}
  \|Du\|_{L^2(B_{1/2})}\leq C \|u\|_{L^1(B_{1})}.
\end{equation*}

Now we consider general $\pmb{b}\in L^p(B_1)$. Take $\pmb{b}^k \in C^{\infty}(\bar{B}_1)$ with $\mathrm{div}~\pmb{b}^k = 0$ such that $\pmb{b}^k\rightarrow \pmb{b}$ in $L^2(B_1)$. Let $v^k$ be the solution of
\begin{equation*}
  \left\{
  \begin{aligned}
    &-\Delta v^k+\pmb{b}^k\cdot Dv^k=0&&~~\mbox{in}~~B_1;\\
    &v^k=u&&~~\mbox{on}~~\partial B_1.
  \end{aligned}
  \right.
\end{equation*}
Then $v^k \in C^{\infty}(B_1)\cap H^1 (\bar{B}_1)$ and from above argument,
\begin{equation}\label{e.vk2}
  \|D v^k\|_{L^2(B_{1/2})}\leq C\|v^k\|_{L^1(B_1)}.
\end{equation}

Let $w^k=v^k-u$. Then $w$ satisfies
\begin{equation*}
  \left\{
  \begin{aligned}
    &-\Delta w^k+\pmb{b}^k\cdot Dw^k=(\pmb{b}-\pmb{b}^k)\cdot Du&&~~\mbox{in}~~B_1\\
    &w^k=0&&~~\mbox{on}~~\partial B_1.
  \end{aligned}
  \right.
\end{equation*}
By \cite[Proposition 2.4]{MR2031029} (see also \cite[Theorem 3.1]{MR3818670}), we have
\begin{equation*}
  \|w^k\|_{L^1(B_1)}\leq C\|\pmb{b}-\pmb{b}^k\|_{L^2(B_1)}\|Du\|_{L^2(B_1)}\rightarrow 0.
\end{equation*}
Hence,
\begin{equation*}
  v^k\rightarrow u ~~\mbox{in}~~L^1(B_1).
\end{equation*}
By \cref{e.vk2}, there exists a subsequence (denoted by $v^k$ again) such that
\begin{equation*}
  \begin{aligned}
    Dv^k\rightharpoonup Du~~\mbox{in}~~L^2(B_{1/2}).
  \end{aligned}
\end{equation*}
Then
\begin{equation*}
  \|Du\|_{L^2(B_{1/2})}\leq \varliminf_{k\to \infty} \|Dv^k\|_{L^2(B_{1/2})}
   \leq C \varliminf_{k\to \infty} \|v^k\|_{L^1(B_{1})}= C\|u\|_{L^1(B_{1})}\leq C \|u\|_{L^2(B_{1})}.
\end{equation*}
~\\

\textbf{Case (2):} $\pmb{b}\in L^2(B_1)$ and
\begin{equation}\label{e2.3}
\int_{B_1} |\pmb{b}|^2 \varphi^2 \leq C_{E} \int_{B_1} |D \varphi|^2,~\forall~\varphi\in C^{\infty}_c (B_1).
\end{equation}
The proof is similar to \textbf{Case (1)} except \Cref{e.1.b}. Since we have assumed $\pmb{b}\in C^{\infty}(\bar{B}_1)$, \cref{e2.3} holds for $\varphi\in H^1_0(B_1)$. We treat \cref{e.1.b} as follows. By \Cref{e2.3} and \Cref{e2.11},
\begin{equation*}
  \begin{aligned}
    \left|\int_{B_1}b^i u^2 \psi \psi_i\right|
    &\leq \left(\int_{B_1} |\pmb{b}|^2 u^2 \psi^2\right)^{1/2}
    \left(\int_{B_1}u^2|D\psi|^2\right)^{1/2}\\
    &\leq \|u D\psi\|_{L^2(B_1)}\left(C_E \int_{B_1}|\psi Du|^2+|u D\psi|^2\right)^{1/2}\\
    &\leq C\|u D\psi\|_{L^2(B_1)}\left(\|\psi Du\|_{L^2(B_1)}+\|u D\psi\|_{L^2(B_1)} \right)\\
    &\leq \frac{1}{8}\|\psi Du\|^2_{L^2(B_1)}+\bar{C} \|u D\psi\|^2_{L^2(B_1)}\\
    &\leq \frac{1}{4}\|\psi Du\|^2_{L^2(B_1)}+C \|u\|^2_{L^1(B_1)}.
  \end{aligned}
\end{equation*}
~\\

\textbf{Case (3):} $\pmb{b}\in BMO^{-1}(B_1)$. Take $\psi^4 u$ as test function in \cref{e.weak2.1} to get
\begin{equation*}
  \int_{B_1} u_i (4\psi^3 \psi_i u+\psi^4 u_i)-(A^{ij}-A_{B_1}^{ij})u_i(\psi^4u)_j=0,
\end{equation*}
where $A_{B_1}^{ij}=\fint_{B_1} A^{ij}$. That is,
\begin{equation}\label{e.energy.4}
  \int_{B_1} u_i (4\psi^3 \psi_i u+\psi^4 u_i)
  -\frac{1}{2}(A^{ij}-A_{B_1}^{ij})\left(u_i(\psi^4u)_j-u_j(\psi^4u)_i\right)=0.
\end{equation}
Note that
\begin{equation}\label{e1.7}
  \begin{aligned}
\int_{B_1} &(A^{ij}-A_{B_1}^{ij})\left(u_i(\psi^4u)_j-u_j(\psi^4u)_i\right)\\
&=\int_{B_1} (A^{ij}-A_{B_1}^{ij})\left(4\psi^3\psi_j u u_i+ \psi^4 u_i u_j-4\psi^3 \psi_i u u_j-\psi^4u_iu_j\right)\\
&=4\int_{B_1} (A^{ij}-A_{B_1}^{ij})\left(\psi^3\psi_j u u_i-\psi^3 \psi_i u u_j\right).
  \end{aligned}
\end{equation}
Since $\psi$ is smooth, the last term can be understood in the usual integral sense. By combining with $A^{ij}=-A^{ji}$, we have
\begin{equation*}
  \int_{B_1} u_i (4\psi^3 \psi_i u+\psi^4 u_i)-4(A^{ij}-A_{B_1}^{ij})\psi^3\psi_j u_i u=0.
\end{equation*}
From the Young's inequality,
\begin{equation*}
  \left|\int_{B_1} \psi^3 \psi_i u u_i \right|\leq \frac{1}{8} \int_{B_1} \psi^4 u_i^2+C \int_{B_1} u^2
  =\frac{1}{8}\|\psi^2 Du\|_{L^{2}(B_1)}^2+C\|u\|_{L^2(B_1)}^2.
\end{equation*}
By the H\"{o}lder inequality,
\begin{equation*}\label{e2.1}
  \left|\int_{B_1} (A^{ij}-\bar{A}^{ij})\psi^3 \psi_j u u_i \right|
  \leq \left(\int_{B_1}|A^{ij}-\bar{A}^{ij}|^q\right)^{\frac{1}{q}}
  \left(\int_{B_1} \psi^4 u_i^2 \right)^{\frac{1}{2}}
  \left(\int_{B_1}|\psi \psi_j u|^{q_0} \right)^{\frac{1}{q_0}},
\end{equation*}
where $1/q+1/q_0=1/2$. Since $q>n$, we can choose $0<\delta_0<1/2$ such that
\begin{equation*}
  2^*=\frac{2n}{n-2}=q_0(1-\delta_0)\cdot \frac{2}{2-q_0\delta_0}.
\end{equation*}
Next, we take $\psi$ as in \textbf{Case (1)} such that
\begin{equation*}
\frac{|D\psi|}{\psi^{1-2\delta_0}}\leq C.
\end{equation*}
Then by the H\"{o}lder inequality and Sobolev inequality,
\begin{equation*}
  \begin{aligned}
   \left(\int_{B_1}|\psi \psi_j u|^{q_0} \right)^{\frac{1}{q_0}}
   &=\left(\int_{B_1}(\psi^{2(1-\delta_0)} |u|^{1-\delta_0}
   \frac{\psi_j}{\psi^{1-2\delta_0}} |u|^{\delta_0})^{q_0} \right)^{\frac{1}{q_0}}\\
   &\leq C\|u\|_{L^2(B_1)}^{\delta_0}\|\psi^2 u\|_{L^{2^*}(B_1)}^{1-\delta_0}\\
   &\leq C\left(\|u\|_{L^2(B_1)}^{\delta_0}\|\psi^2 Du\|_{L^{2}(B_1)}^{1-\delta_0}+\|u\|_{L^2(B_1)}\right).
  \end{aligned}
\end{equation*}
Hence,
\begin{equation*}
  \left|\int_{B_1} (A^{ij}-\bar{A}^{ij})\psi^3 \psi_j u u_i \right|
  \leq \frac{1}{8}\|\psi^2 Du\|_{L^{2}(B_1)}^{2}+C\|u\|_{L^2(B_1)}^2,
\end{equation*}
where $C$ depends only on $n$ and $\|\pmb{A}\|_{BMO(0)}$.
By combining the above estimates, we have
\begin{equation*}
\|Du\|_{L^{2}(B_{1/2})}\leq C\|u\|_{L^2(B_1)}
\end{equation*}
where $C$ depends only on $n$ and $\|\pmb{b}\|_{BMO^{-1}(0)}$.~\qed~\\

\begin{remark}\label{e.energy.2}
If $\pmb{b}$ is smooth enough (e.g. $\pmb{b}\in L^{\infty}(B_1)$), the energy inequality holds with $\mathrm{div}~\pmb{b} \leq 0$ since we can take $\varphi^2 u$ as the test function directly.
\end{remark}

\begin{remark}\label{e.energy.0}
Note that the constant $C$ in \cref{e1.6} depends only on $\|\pmb{b}\|_{BMO^{-1}(0)}$ rather than $\|\pmb{b}\|_{BMO^{-1}(B_1)}$. The latter is assumed in \Cref{ener-e} and it is used to guarantee that the definition of the weak solution makes sense (see \cref{e.weak2}). However, the troublesome term $A^{ij}u_i u_j$ is canceled in the calculation (see \cref{e1.7}). Hence, the constant $C$ doesn't depend on $\|\pmb{b}\|_{BMO^{-1}(B_1)}$ ultimately.
\end{remark}


\begin{remark}\label{re.energy.1}
With (3), Li and Pipher \cite[Corollary 3.1]{MR3936348} gave a energy inequality for the elliptic equation with the righthand $\|u\|_{L^p}$ for $p>2$ instead of $\|u\|_{L^2}$.
\end{remark}

\begin{remark}\label{re.2.5}
We have used the following fact in the proof: by choosing a proper $\psi$ such that $|D \psi|/\psi^{1-\delta_0} \leq C$. This technique is inspired by \cite[P. 95-96]{MR2760150} and \cite[P. 254-255]{MR2031029}.
\end{remark}
~\\

In the following, we prove the pointwise $C^{\alpha}$ regularity for \cref{e.equation}. Here we just prove the result with condition \cref{e1.2} because the proofs are similar for the conditions \cref{e1.3} and \cref{e1.4}. We point out the differences of proofs in remarks.

First, we prove the key step towards interior $C^{\alpha}$ regularity.
\begin{lemma}\label{l-Ca-key}
For any $0<\alpha<1$, there exists $\delta>0$ depending only on $n,p$ and $\alpha$ such that if $u$ is a weak solution of \Cref{e.equation} with
\begin{equation}\label{e.small}
  \begin{aligned}
    &\|u\|^*_{L^{2}(B_1)}\leq 1,~~\|\pmb{b}\|^*_{L^{p}(B_1)}\leq \delta,
  \end{aligned}
\end{equation}
then there exists a constant $\bar P$ such that
\begin{equation}\label{e-lCa-udis}
  \|u-\bar P\|^*_{L^2(B_{\eta})}\leq \eta^{\alpha}
\end{equation}
and
\begin{equation}\label{e-lca-p}
|\bar P|\leq \bar C,
\end{equation}
where $\bar C$ and  depends only on $n$, and $0<\eta<1/2$ depends also on $\alpha$.
\end{lemma}

\proof We prove the lemma by contradiction. Suppose that the lemma is false. Then there exist $0<\alpha<1$ and sequences of $u_m$ and $\pmb{b}_m$ such that
\begin{equation}\label{e.3.um}
\Delta u_m+b_m^iu_{m,i}=0~~\mbox{in}~~B_1
\end{equation}
with
\begin{equation*}
  \begin{aligned}
    &\|u_m\|^*_{L^{2}(B_1)}\leq 1,~~
    \|\pmb{b}_m\|^*_{L^{p}(B_1)}\leq 1/m.
  \end{aligned}
\end{equation*}
In addition, for any constant $P$ with $|P|\leq \bar C$,
\begin{equation}\label{e-lCa-1}
  \|u_m-P\|^*_{L^{2}(B_{\eta})}> \eta^{\alpha},
\end{equation}
where $\bar C$ is to be specified later and $0<\eta<1/2$ is taken small such that
\begin{equation}\label{e-lCa-2}
\bar C\eta^{1-\alpha}<1/2.
\end{equation}

By $\|u_m\|^*_{L^2(B_1)}\leq 1$ and the energy inequality (\Cref{ener-e}), ${u_m}$ is bounded in $H^1(B_{3/4})$. Then there exists a subsequence (denoted by $u_m$ again) such that $u_m$ converges strongly in $L^2$ and weekly in $H^1$ to some function $\bar{u}$ in $B_{3/4}$. Then $\|\bar{u}\|^*_{L^2(B_{3/4})}\leq (4/3)^n$. Moreover, since $\|\pmb{b}_m\|^*_{L^{p}(B_1)}\to 0$, $\bar{u}$ is a weak solution of
\begin{equation}\label{e.laplace}
  \Delta \bar{u}=0~~\mbox{in}~~B_{3/4}.
\end{equation}
Then $\bar u$ is smooth. By the $C^1$ estimate for \Cref{e.laplace}, there exists a constant $\bar{P}$ such that
\begin{equation*}
  \|\bar u-\bar{P}\|^*_{L^2(B_r)}\leq  \bar Cr, ~~\forall ~0<r<1/2
\end{equation*}
and
\begin{equation*}
  |\bar{P}|\leq \bar C,
\end{equation*}
where $\bar{C}$ depends only on $n$.

By noting \cref{e-lCa-2}, we have
\begin{equation}\label{e-lCa-3}
  \|u-\bar{P}\|^*_{L^{2}(B_{\eta})}\leq \eta^{\alpha}/2.
\end{equation}
However, from \cref{e-lCa-1},
\begin{equation*}
  \|u_m-\bar{P}\|^*_{L^{2}(B_{\eta})}> \eta^{\alpha}.
\end{equation*}
Let $m\rightarrow \infty$, we have
\begin{equation*}
    \|u-\bar{P}\|^*_{L^{2}(B_{\eta})}\geq \eta^{\alpha},
\end{equation*}
which contradicts with \cref{e-lCa-3}.  ~\qed~\\

\begin{remark}\label{re.3.1}
If we consider the condition \cref{e1.3}, we assume
\begin{equation*}
  C_E\leq \delta
\end{equation*}
in above lemma. That is,
\begin{equation}\label{e2.2}
  \int_{B_1} |\pmb{b}|^2\varphi^2 \leq \delta \int_{B_1} |D\varphi|^2,~~\forall~~\varphi \in C^{\infty}_c(B_1).
\end{equation}

Similarly, we can prove the lemma by contradiction. Then there exist $0<\alpha<1$ and sequences of $u_m$, $\pmb{b}_m$ such that
\begin{equation*}
\|u_m\|^*_{L^{2}(B_1)}\leq 1,~\int_{B_1} |\pmb{b}_m|^2\varphi^2 \leq \frac{1}{m}\int_{B_1} |D\varphi|^2
,~~\forall~~\varphi \in C^{\infty}_c(B_1)
\end{equation*}
and \crefrange{e.3.um}{e-lCa-1} hold.

By the compactness (energy inequality), $u_m$ converges to some $\bar{u}$ in $B_{3/4}$. Then we can verify that $\bar{u}$ is a harmonic function in the weak sense as follows. Take $\psi^2$ ($\psi \in C_c^{\infty}(B_{3/4})$) as the test function, we have
\begin{equation*}
  \int_{B_{3/4}} u_{m,i} (\psi^2)_i-b^i_mu_{m,i}\psi^2=0.
\end{equation*}
Note that
\begin{equation*}
  \begin{aligned}
    \left|\int_{B_{3/4}}b^i_mu_{m,i}\psi^2\right|
    &=2\left|\int_{B_{3/4}}b^i_mu_{m}\psi \psi_i\right|
    \leq 2\left(\int_{B_{3/4}}|\pmb{b}_m|^2\psi^2\right)^{\frac{1}{2}}\left(\int_{B_{3/4}}u_{m}^2|D\psi|^2\right)^{\frac{1}{2}}\\
    &\leq C\left(\frac{1}{m} \int_{B_{3/4}}|D\psi|^2\right)^{\frac{1}{2}}\to 0.
  \end{aligned}
\end{equation*}
Then $\bar{u}$ is a weak solution of \Cref{e.laplace}. The rest of the proof is similar to that of \Cref{l-Ca-key} and we omit it.

If we consider the condition \cref{e1.4}, we assume
\begin{equation*}
  \|\pmb{b}\|_{BMO^{-1}(0)}\leq \delta
\end{equation*}
in above lemma. That is,
\begin{equation*}
  \|\pmb{A}\|_{BMO(0)}\leq \delta.
\end{equation*}

We can prove the lemma similarly. We just present here how the $\bar{u}$ satisfies the Laplace's equation. Take $\psi \in C_c^{\infty}(B_{3/4})$ as the test function, we have
\begin{equation*}
  \int_{B_{3/4}} u_{m,i} \psi_i-(A^{ij}_m-A^{ij}_{m,B_1}) u_{m,i}\psi_j=0
\end{equation*}
and
\begin{equation*}
  \|\pmb{A}_m\|_{BMO(0)}\leq 1/m.
\end{equation*}
Note that
\begin{equation*}
  \begin{aligned}
    &\left|\int_{B_{3/4}} (A^{ij}_m-A^{ij}_{m,B_1}) u_{m,i}\psi_j\right|
    \leq C\|\pmb{A}_m\|_{BMO(0)}\|\nabla u_m\|_{L^2(B_{3/4})}\|\nabla \psi\|_{L^{\infty}(B_{3/4})}.
  \end{aligned}
\end{equation*}
Then, by letting $m\rightarrow \infty$, $\bar{u}$ is a weak solution of \Cref{e.laplace}.
\end{remark}
~\\

Now, we give the~\\
\noindent\textbf{Proof of \Cref{t-Ca-elli}.} Without loss of generality, we assume that $\|u\|^*_{L^{2}(B_1)}\leq 1$. To prove that $u$ is $C^{\alpha}$ at $0$, we only need to prove the following. There exist a sequence of constants $P_m$ ($m\geq -1$) such that for all $m\geq 0$,
\begin{equation}\label{e-tCas-u}
\|u-P_m\|^*_{L^{2}(B_{\eta^{m}})}\leq \eta ^{m\alpha}
\end{equation}
and
\begin{equation}\label{e-tCas-P}
|P_m-P_{m-1}|\leq \bar C\eta^{(m-1)\alpha},
\end{equation}
where $\bar C$ and $\eta$ are as in \Cref{l-Ca-key}. We prove \cref{e-tCas-u} and \cref{e-tCas-P} by induction. For $m=0$, by setting $P_0=P_{-1}\equiv 0$, \crefrange{e-tCas-u}{e-tCas-P} hold clearly. Suppose that the conclusion holds for $m\leq m_0$. We need to prove the conclusion for $m=m_0+1$.

Let $r=\eta ^{m_0}$, $y=x/r$ and
\begin{equation}\label{e-tCas-v}
  v(y)=\frac{u(x)-P_{m_0}}{r^{\alpha}}.
\end{equation}
Then $v$ satisfies
\begin{equation}\label{e-tCas-f}
\Delta v+\tilde b^iv_i=0~~\mbox{in}~~B_1,
\end{equation}
where
\begin{equation}\label{e-tCas-new}
  \begin{aligned}
&\tilde {\pmb{b}}(y)=r \pmb{b}(x).
  \end{aligned}
\end{equation}

In the following, we show that \cref{e-tCas-f} satisfies the assumptions of \Cref{l-Ca-key}. First, by \cref{e-tCas-u} and \cref{e-tCas-v},
\begin{equation*}
\|v\|^*_{L^{2}(B_1)}= r^{-\alpha}\|u-P_{m_0}\|^*_{L^{2}(B_r)}\leq 1.
\end{equation*}
Next, by $\|\pmb{b}\|_{C^{-1}(0)}\leq \delta$ and \cref{e-tCas-new},
\begin{equation*}
  \begin{aligned}
   \|\tilde {\pmb{b}}\|^*_{L^{p}(B_1)}=r \|\pmb{b}\|^*_{L^{p}(B_r)}
    \leq r \|\pmb{b}\|_{C^{-1}(0)} r^{-1} \leq \delta.
  \end{aligned}
\end{equation*}

Therefore, \cref{e-tCas-f} satisfies the assumptions of \Cref{l-Ca-key} and hence there exists a constant $\bar P$ such that
\begin{equation*}
\begin{aligned}
    \|v-\bar P\|^*_{L^{2}(B_{\eta})}&\leq \eta ^{\alpha}
\end{aligned}
\end{equation*}
and
\begin{equation*}
|\bar P|\leq \bar C.
\end{equation*}
Let $P_{m_0+1}=P_{m_0}+r^{\alpha}\bar P$. Then \cref{e-tCas-P} hold for $m_0+1$. Recalling \cref{e-tCas-v}, we have
\begin{equation*}
  \begin{aligned}
\|u-P_{m_0+1}\|^*_{L^{2}(B_{\eta^{m_0+1}})}
&= \|u-P_{m_0}-r^{\alpha}\bar P\|^*_{L^{2}(B_{\eta r})}\\
&= \|r^{\alpha}v-r^{\alpha}\bar P\|^*_{L^{2}(B_{\eta})}\\
&\leq r^{\alpha}\eta^{\alpha}=\eta^{(m_0+1)\alpha}.
  \end{aligned}
\end{equation*}
Hence, \cref{e-tCas-u} holds for $m=m_0+1$. By induction, the proof is completed.\qed~\\

\begin{remark}\label{re.3.3}
If we consider the condition \cref{e1.3}, the proof can be proceeded in a similar way. The main difference is to verify that $\tilde{\pmb{b}}$ satisfies \cref{e2.2}. Indeed, for any $\varphi(y) \in C^{\infty}_c(B_1)$,
\begin{equation*}
  \begin{aligned}
    &\int_{B_1} |\tilde{\pmb{b}}|^2\varphi^2 dy=\int_{B_r} |r\pmb{b}|^2\tilde{\varphi}^2 \frac{1}{r^n}dx
    \leq \delta \int_{B_r} r^2|D\tilde{\varphi}|^2 \frac{1}{r^n}dx=\delta \int_{B_1} |D\varphi|^2 dy,
  \end{aligned}
\end{equation*}
where $\tilde{\varphi}(x)=\varphi(y)$.

Similarly, for the condition \cref{e1.4}, we have $\tilde{\pmb{A}}(y)=\pmb{A}(x)$ during the scaling and
\begin{equation*}
  \begin{aligned}
    &\|\tilde{\pmb{A}}\|_{BMO(0)}
    =\sup_{0<\rho<1} \left(\fint_{B_{\rho}} |\tilde{A}^{ij}-\tilde{A}^{ij}_{B_{\rho}}|^q \right)^{\frac{1}{q}}\\
    &=\sup_{0<\rho <1} \left(\fint_{B_{\rho r}} |A^{ij}-A^{ij}_{B_{\rho r}}|^q \right)^{\frac{1}{q}}
    \leq \|\pmb{A}\|_{BMO(0)}\leq \delta.
  \end{aligned}
\end{equation*}
That is, $\|\tilde{\pmb{b}}\|_{BMO^{-1}(0)}\leq \delta$, which meets the requirement in \Cref{l-Ca-key}.
\end{remark}
~\\

\section{$C^{\alpha}$ regularity for parabolic equations}\label{sec:2}
In this section, we prove the $C^{\alpha}$ regularity for the parabolic case. Since the proofs are similar to that for elliptic equations, we only give the outline of the proof and omit the details. In this section, if we say that $u$ is a weak solution of \Cref{e.equation}, it means that the weak solution of (\textbf{P}) in \Cref{e.equation}. Similar to the elliptic equation, we first give the energy inequality for parabolic equations. The proof is similar to that of \Cref{ener-e} and we omit the details here.

\begin{lemma}\label{ener-p}
  Let $u$ be a weak solution of \Cref{e.equation}. Assume that $\pmb{b}$ satisfies one of the following conditions:
  \begin{equation}\label{c.b.PE}
    \begin{aligned}
      &(1)~\pmb{b}\in L^{p,p'}(Q_1)~~(1<\frac{n}{p}+\frac{2}{p'}<2);\\
      &(2)~\pmb{b}\in L^2(Q_1),~~\int_{Q_1} |\pmb{b}|^2 \varphi^2 \leq C_{P} \int_{Q_1} |D \varphi|^2,
      ~~\forall~\varphi \in C^{\infty}_c (Q_1);\\
      &(3)~\pmb{b}\in L^{\infty}(-1,0;BMO^{-1}(B_1)).
    \end{aligned}
  \end{equation}
  Then
  \begin{equation*}
    \|u\|_{L^2(-1/4,0;H^1(B_{1/2}))}+\|u\|_{L^{\infty}(-1/4,0;L^2(B_{1/2}))}\leq C \|u\|_{L^2(Q_{1})},
  \end{equation*}
  where $C$ depends only on $n$ and $\|\pmb{b}\|_{L^{p,p'}(Q_1)}$ (or $C_P$ or $\|\pmb{b}\|_{L^{\infty}(-1,0;BMO^{-1}(0))}$).
\end{lemma}

Since $H^{1}(B_1)$ is compactly embedded into $L^2(B_1)$, the energy inequality provides the compactness for elliptic equations. For the parabolic equation, the space $L^2(-1,0;H^1(B_1))\cap L^{\infty}(-1,0;L^2(B_1))$ is not embedded compactly into $L^2(Q_1)$. Hence, not only do we need the energy inequality, but also we should make estimates of $u_t$ for the compactness.

Before estimating $u_t$, let's introduce
\begin{definition}\label{df-weak-u}
Let $X$ be a Banach space and $u,v\in L^1(T_1,T_2;X)$. We say that $v$ is the weak derivative of $u$ with respect to $t$ and denote $v=u_t$ if
\begin{equation}\label{e.weaku.1}
\int_{T_1}^{T_2} u(t)\eta '(t)dt=-\int_{T_1}^{T_2} v(t)\eta(t)dt,~\forall ~\eta\in C^{\infty}_c(T_1,T_2;\mathbb{R}).
\end{equation}
\end{definition}

Next, the estimate of $u_t$ is
\begin{lemma}\label{th.ut.1}
Let $u\in L^2(-1,0;H^1(B_1))\cap L^{\infty}(-1,0;L^2(B_1))$ be a weak solution of \Cref{e.equation} and $\pmb{b}$ satisfy one of the conditions in \Cref{c.b.PE}. Then $u_t\in L^{\tilde{p}}(-1/4,0;H^{-1}(B_{1/2}))$ for
\begin{equation*}
\tilde{p}=\left\{
  \begin{aligned}
    &2(n/p+2/p')^{-1},&&~~\mbox{if (1) holds};\\
    &2,&&~~\mbox{if (2) or (3) holds }~~
  \end{aligned}
\right.
\end{equation*}
and
%
\begin{equation*}
\|u_t\|_{L^{\tilde{p}}(-1/4,0;H^{-1}(B_{1/2}))}\leq C\|u\|_{L^2(Q_1)},
\end{equation*}
where $C$ depends also on $n$ and $\|\pmb{b}\|_{L^{p,p'}(Q_1)}$ (or $C_P$, or $\|\pmb{b}\|_{L^{\infty}(-1,0;BMO^{-1}(0))}$).
\end{lemma}

\proof We prove the lemma according to each case in \Cref{c.b.PE}. \textbf{Case (1):} We first show that
\begin{equation}\label{e2.9}
u_t=\Delta u-b^iu_i:=v
\end{equation}
in the weak sense of \Cref{df-weak-u}. Indeed, for any $\eta\in C^{\infty}_c(-1,0)$ and $\varphi\in C^{\infty}_c(B_1)$, take $\eta\varphi$ as the test function. Then by the definition of weak solution,
\begin{equation*}
\begin{aligned}
\int_{Q_1}u\eta_t\varphi=&\int_{Q_1}u_i \eta\varphi_i+b^iu_i\eta\varphi.
\end{aligned}
\end{equation*}
That is,
\begin{equation*}
\int_{-1}^{0} \langle u\eta_t,\varphi \rangle=\int_{-1}^{0} \langle v\eta,\varphi \rangle,
\end{equation*}
where $\langle \cdot,\cdot\rangle$ denotes the pairing between $H^{-1}(B_1)$ and $H^1_0(B_1)$. Then
\begin{equation*}
 \langle \int_{-1}^{0}u\eta_t,\varphi \rangle= \langle \int_{-1}^{0}v\eta,\varphi \rangle,
\end{equation*}
which implies
\begin{equation*}
 \int_{-1}^{0}u\eta_t= \int_{-1}^{0}v\eta ~~\mbox{ in}~H^{-1}(B_1).
\end{equation*}
Thus, we have proved \cref{e2.9}.

Next, we show the estimate for $u_t$. For any $\varphi \in C_c^{\infty}(B_{1/2})$,
\begin{equation}\label{e.ut.e}
  \begin{aligned}
\big|\langle u_t,\varphi\rangle\big|
=&\bigg|\int_{B_{1/2}}u_i\varphi_i+b^iu_i\varphi\bigg|
= \bigg|\int_{B_{1/2}}u_i\varphi_i-b^iu \varphi_i\bigg|\\
\leq& C\left(\|Du\|_{L^2(B_{1/2})}\|D\varphi\|_{L^2(B_{1/2})}
+\|\pmb{b}\|_{L^p(B_{1/2})}\|u\|_{L^{p_0}(B_{1/2})}\|D \varphi\|_{L^2(B_{1/2})}\right)\\
\leq &C\left(\|Du\|_{L^2(B_{1/2})}+\|\pmb{b}\|_{L^p(B_{1/2})}\|u\|_{L^{p_0}(B_{1/2})}\right)\|D \varphi\|_{L^2(B_{1/2})},
  \end{aligned}
\end{equation}
where $C$ depends only on $n$ and
\begin{equation}\label{e.p.p0}
\frac{1}{p}+\frac{1}{p_0}=\frac{1}{2}.
\end{equation}
Thus,
\begin{equation*}
  \|u_t(\cdot,t)\|_{H^{-1}(B_{1/2})}\leq C\left(\|Du\|_{L^2(B_{1/2})}+\|\pmb{b}\|_{L^p(B_{1/2})}\|u\|_{L^{p_0}(B_{1/2})}\right).
\end{equation*}
By noting the definition of $\tilde{p}$, i.e.,
\begin{equation*}
  \frac{n}{p}+\frac{2}{p'}=\frac{2}{\tilde{p}}
\end{equation*}
and \Cref{e.p.p0}, we have
\begin{equation*}
  1<\tilde{p}<2~~\mbox{and}~~\frac{1}{p'}+\frac{1}{p'_0}=\frac{1}{\tilde{p}}.
\end{equation*}
Next, by the H\"{o}lder inequality and the G-N inequality \cite{MR4513000},
\begin{equation}\label{e4.1}
\begin{aligned}
  \|u_t\|_{L^{\tilde{p}}(-1/4,0;H^{-1}(B_{1/2}))}
  \leq& C\left(\|u\|_{L^2(-1/4,0;H^1(B_{1/2}))}
   +\|\pmb{b}\|_{L^{p,p'}(Q_{1/2})}\|u\|_{L^{p_0,p_0'}(-1/4,0;H^1(B_{1/2}))}\right)\\
  \leq&C\left(\|u\|_{L^2(-1/4,0;H^1(B_{1/2})}+\|u\|_{L^{\infty}(-1/4,0;L^2(B_{1/2}))}\right),
\end{aligned}
\end{equation}
where $C$ depends only on $n$ and $\|\pmb{b}\|_{L^{p,p'}(Q_{1/2})}$.

Finally, by combining with the energy inequality (\Cref{ener-p}),
\begin{equation*}
\begin{aligned}
  \|u_t\|_{L^{\tilde{p}}(-1/4,0;H^{-1}(B_{1/2}))}\leq C\|u\|_{L^2(Q_1)}.
\end{aligned}
\end{equation*}
~\\

\textbf{Case (2):} The proof of is similar to \textbf{Case (1)} and we only show the estimate of $u_t$. For any $\varphi \in L^2(-1/4,0;H_0^1(B_{1/2}))$,
 \begin{equation*}
  \begin{aligned}
\bigg|\int_{-1/4}^{0} \langle u_t,\varphi \rangle \bigg|
=&\bigg|\int_{-1/4}^{0} \langle \Delta u-b^iu_i,\varphi \rangle \bigg|
= \bigg|\int_{Q_{1/2}}u_i\varphi_i+b^iu_i\varphi\bigg|\\
\leq& \bigg|\int_{Q_{1/2}}u_i\varphi_i\bigg|+\left(\int_{Q_{1/2}}|\pmb{b}|^2 \varphi^2 \right)^{\frac{1}{2}} \left(\int_{Q_{1/2}} u_i^2\right)^{\frac{1}{2}}\\
\leq& \|Du\|_{L^2(Q_{1/2})}\|D\varphi\|_{L^2(Q_{1/2})}
+C_P^{\frac{1}{2}}\|D \varphi\|_{L^2(Q_{1/2})}\|Du\|_{L^{2}(Q_{1/2})}\\
\leq &(1+C_P^{\frac{1}{2}})\|Du\|_{L^2(Q_{1/2})}\|D \varphi\|_{L^2(Q_{1/2})}.
  \end{aligned}
\end{equation*}
Then with the aid of the energy inequality,
\begin{equation*}
  \|u_t\|_{L^2(-1/4,0;H^{-1}(B_{1/2}))}\leq C\|u\|_{L^2(Q_1)}.
\end{equation*}
~\\

\textbf{Case (3):} Similar to \textbf{Case (2)}, we just show the estimate of $u_t$. Let $\varphi$ be as in \textbf{Case (1)}.
\begin{equation*}
  \begin{aligned}
\big|\langle u_t,\varphi\rangle\big|
=&\bigg|\int_{B_{1/2}}u_i\varphi_i-\frac{1}{2}A^{ij}(u_i\varphi_j-u_j\varphi_i)\bigg|\\
\leq& \|Du\|_{L^2(B_{1/2})}\|D\varphi\|_{L^2(B_{1/2})}
+C\|\pmb{A}\|_{BMO(B_{1/2})}\|D u\|_{L^{2}(B_{1/2})}\|D \varphi\|_{L^2(B_{1/2})}\\
\leq& C\left(\|Du\|_{L^2(B_{1/2})}
+\|\pmb{A}\|_{BMO(B_{1/2})}\|D u\|_{L^{2}(B_{1/2})}\right)\|D\varphi\|_{L^2(B_{1/2})},
  \end{aligned}
\end{equation*}
where $C$ depends only on $n$.
%
~\qed~\\

Once the estimate of $u_t$ is built, we have the necessary compactness for weak solutions. Let's recall the Aubin-Lions Lemma (see \cite[Theorem 1.3]{MR3101774}).
\begin{lemma}\label{le3.1}
Suppose that $X_0,X$ and $X_1$ are Banach spaces satisfying $X_0\subset X \subset X_1$ and $X_0\hookrightarrow X$ compactly. Let $T>0$, $\alpha_0,\alpha_1\in (1,+\infty)$ and
\begin{equation*}
\mathcal{Y}=\left\{v\in L^{\alpha_0}(-T,0;X_0): v_t\in L^{\alpha_1}(-T,0;X_1)\right\}.
\end{equation*}
Note that $\mathcal{Y}$ is a Banach space with the norm
\begin{equation*}
  \|v\|_{\mathcal{Y}}:=\|v\|_{L^{\alpha_0}(-T,0;X_0)}+\|v_t\|_{L^{\alpha_1}(-T,0;X_1)}.
\end{equation*}

Then the embedding $\mathcal{Y}\hookrightarrow L^{\alpha_0}(-T,0;X_0)$ is compact.
\end{lemma}

Now, we apply \Cref{le3.1} with $X_0=H^1(B_{1/2}),X=L^2(B_{1/2}), X_1=H^{-1}(B_{1/2})$, $T=1/4,\alpha_0=2,\alpha_1=\tilde{p}$ (note that $\tilde{p}>1$). By the energy inequality (\Cref{ener-p}) and the estimate of $u_t$ (\Cref{th.ut.1}), we know that any sequence of weak solutions $u_m$ with a uniform $L^2$ norm bound (i.e., $\|u_m\|_{L^2(Q_1)}\leq C,~\forall ~m$) is precompact in $L^2(Q_{1/2})$. Then we can prove the following key step in an exactly the same way as the elliptic equation.
\begin{lemma}\label{l.key.p}
For any $0<\alpha<1$, there exists $\delta>0$ depending only on $n,p$ and $\alpha$ such that if $u$ is a weak solution of \cref{e.equation} with
\begin{equation*}
  \begin{aligned}
    &\|u\|^*_{L^{2}(Q_1)}\leq 1,~~\|\pmb{b}\|^*_{L^{p,p'}(Q_1)}\leq \delta~~(\mbox{or}~~C_P\leq \delta
    ~~\mbox{or}~~\|\pmb{b}\|_{L^{\infty}(-1,0;BMO^{-1}(0))}\leq \delta).
  \end{aligned}
\end{equation*}
Then there exists a constant $\bar P$ such that
\begin{equation}\label{e-lCa-udis}
  \|u-\bar P\|^*_{L^2(Q_{\eta})}\leq \eta^{\alpha}
\end{equation}
and
\begin{equation}\label{e-lca-p}
|\bar P|\leq \bar C,
\end{equation}
where $\bar C$ and  depends only on $n$, and $0<\eta<1/2$ depends also on $\alpha$.
\end{lemma}

Finally, we can prove \Cref{t-Ca-para} by a scaling argument analogous to the elliptic case. We omit the details.
~\\

\bibliographystyle{amsplain}
\bibliography{PDE2}

\end{document}